\documentclass[11pt]{article}

\usepackage{latexsym}
\usepackage{bezier}

\begin{document}

%%%%%%%%%%%%%%%%%%%%%%%%%%%%%%%%%%%%%%%%%%%%%%%%%%%%%%%%%%%% begin lama.sty
%  LAMA.STY  (auszug)                                           6. M\"arz 92
%

% Spezielle mathematische Symbole ----------------------------------------
\newcommand{\EE}{\mathop{\rm I\! E}\nolimits}
\newcommand{\E}{\mathop{\rm E}\nolimits}
\newcommand{\I}{\mathop{\rm Im\, }\nolimits}
\newcommand{\Str}{\mathop{\rm Str\, }\nolimits}
\newcommand{\Sdet}{\mathop{\rm Sdet\, }\nolimits}
\newcommand{\STr}{\mathop{\rm STr\, }\nolimits}
\newcommand{\R}{\mathop{\rm Re\, }\nolimits}
\newcommand{\CC}{\mathop{\rm C\!\!\! I}\nolimits}
\newcommand{\FF}{\mathop{\rm I\! F}\nolimits}
\newcommand{\KK}{\mathop{\rm I\! K}\nolimits}
\newcommand{\LL}{\mathop{\rm I\! L}\nolimits}
\newcommand{\MM}{\mathop{\rm I\! M}\nolimits}
\newcommand{\NN}{\mathop{\rm I\! N}\nolimits}
\newcommand{\PP}{\mathop{\rm I\! P}\nolimits}
\newcommand{\QQ}{\mathop{\rm I\! Q}\nolimits}
\newcommand{\RR}{\mathop{\rm I\! R}\nolimits}
\newcommand{\ZZ}{\mathop{\sl Z\!\!Z}\nolimits}
% ------------------------------------------------------------------------
\newcommand{\integer}{\mathop{\rm int}\nolimits}
\newcommand{\erf}{\mathop{\rm erf}\nolimits}
\newcommand{\diag}{\mathop{\rm diag}\nolimits}
\newcommand{\fl}{\mathop{\rm fl}\nolimits}
\newcommand{\eps}{\mathop{\rm eps}\nolimits}
\newcommand{\var}{\mathop{\rm var}\nolimits}

%%%%%%%%%%%%%%%%%%%%%%%%%%%%%%%%%%%%%%%%%%%%%%%%%%%%%%%%%%%% end lama.sty

\newcommand{\pfeil}{\rightarrow}

\newcommand{\kat}{{\cal C}}
\newcommand{\rmat}{{\cal R}}
\newcommand{\oalg}{{\cal A}}
\newcommand{\falg}{{\cal F}}
\newcommand{\eich}{{\cal G}}
\newcommand{\hilb}{{\cal H}}
\newcommand{\calm}{{\cal M}}
\newcommand{\mod}{{\cal M}}
\newcommand{\kegel}{{\cal O}}
\newcommand{\kegels}{{\cal K}}
\newcommand{\bigrho}{\rho_\oplus}
\newcommand{\bigphi}{\phi_\oplus}

\newcommand{\tprod}{\otimes}

\newcommand{\horab}{\rule[-1mm]{0pt}{5mm}}
\newcommand{\iso}{\stackrel{\sim}{=}}
\newcommand{\quer}[1]{\overline{#1}}
\newcommand{\schlange}[1]{\widetilde{#1}}
\newcommand{\CVO}[3]{{{#3 \choose #1\hspace{3pt}#2}}}
\newcommand{\clebsch}[6]{{\left[{#1\atop#4}\hspace{3pt}{#2\atop#5}
                                     \hspace{3pt}{#3\atop#6}\right]}}
\newcommand{\sjsymbol}[6]{{\left\{{#1\atop#4}\hspace{3pt}{#2\atop#5}
                                     \hspace{3pt}{#3\atop#6}\right\}}}
\newcommand{\spann}{{{\rm span}}}
\newcommand{\Nat}{{{\rm Nat}}}
\newcommand{\Mor}{{{\rm Mor}}}
\newcommand{\End}{{{\rm End}}}
\newcommand{\Obj}{{{\rm Obj}}}
\newcommand{\ev}{{{\rm ev}}}
\newcommand{\coev}{{{\rm coev}}}
\newcommand{\id}{{{\rm id}}}
\newcommand{\Id}{{{\rm Id}}}
\newcommand{\Vec}{{{\rm Vec}}}
\newcommand{\Rep}{{{\rm Rep}}}
\newcommand{\Amp}{{{\rm Amp}}}
\newcommand{\RAmp}{{{\rm RAmp}}}
\newcommand{\ket}[1]{|#1\rangle}
\newcommand{\vak}{\ket{0}}

\newcommand{\ZA}{{{\rm ZA}}}
\newcommand{\ZB}{{{\rm ZB}}}
\newcommand{\WA}{{{\rm WA}}}
\newcommand{\WB}{{{\rm WB}}}
\newcommand{\TA}{{{\rm TA}}}
\newcommand{\TB}{{{\rm TB}}}

\newcommand{\Ralpha}{{\overline{\alpha}}}
\newcommand{\Rbeta}{{\overline{\beta}}}

\newcommand{\QG}{{{\cal U}_q(sl_2)}}
\newcommand{\QW}{{{\overline \QG}}}

\newenvironment{bew}{Proof:}{\hfill$\Box$}

\newtheorem{bem}{Remark}
\newtheorem{bsp}{Example}
\newtheorem{axiom}{Axiom}
\newtheorem{de}{Definition}
\newtheorem{satz}{Proposition}
\newtheorem{lemma}[satz]{Lemma}
\newtheorem{kor}[satz]{Corollary}
\newtheorem{theo}[satz]{Theorem}

\newcommand{\sbegin}[1]{\small\begin{#1}}
\newcommand{\send}[1]{\end{#1}\normalsize}

\sloppy
%%%%%%%%%%%%%%%%%%%%%%%%%%%%%%%%%%%%%%%%%%%%%%%%%%%%%%%%%%%%%%%%%%%%%%%%%

\title{The Cylinder braiding of the quantum Weyl group of $sl_2$}
\author{Reinhard H\"aring-Oldenburg\\ 
{\small Mathematisches Institut, Bunsenstr. 3-5,
 37073 G\"ottingen, Germany}\\
{\small email: haering@cfgauss.uni-math.gwdg.de}}
\date{August 13, 1996}
\maketitle

\begin{abstract}
Is is shown that the quantum Weyl group of $sl_2$ contains an element
that is a cylinder twist, i.e. it gives rise to representations
of the braid group of Coxeter type B.
\end{abstract}

\section{Introduction}

Every Coxeter graph defines a braid group that is an infinite
covering of its Coxeter group. 	 T. tom Dieck
initiated in \cite{tD1} the systematic study of these braid groups
and their quotient algebras for all root systems.

The Coxeter group $\WA_n$ of type $A_n$  is the permutation group
and the braid group $\ZA_n$ is Artin's braid group. 
For type $B_n$
the Weyl group $\WB_n$ is a semi direct product of the permutation group 
$\WA_n$ with $\ZZ_2^n$.
\begin{de} The braid group $\ZB_n$ of Coxeter type B is generated
by $\tau_0,\tau_1,\ldots,\tau_{n-1}$ with relations
\begin{eqnarray}
\tau_i\tau_j&=&\tau_j\tau_i\quad\mbox{if}\quad |i-j|>1\\
\tau_i\tau_j\tau_i&=&\tau_i\tau_j\tau_i\quad\mbox{if}\quad i,j\geq1,|i-j|=1\\
\tau_0\tau_1\tau_0\tau_1&=&\tau_1\tau_0\tau_1\tau_0	\\
\tau_0\tau_i&=&\tau_i\tau_0\quad i\geq2
\end{eqnarray}
$\tau_0$ is called the cylinder twist. 
\end{de}
Generators $\tau_i,i\geq1$  satisfy the relations of Artin's braid group.

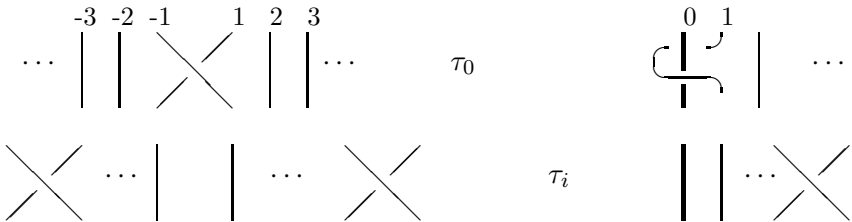
\begin{figure}[ht]
\begin{picture}(150,30)
\put(59,20){\mbox{$\tau_0$}}

\linethickness{0.2mm}
\put(2,20){\mbox{$\cdots$}}
\put(9,26){\mbox{{\small -3}}}
\put(14,26){\mbox{{\small -2}}}
\put(19,26){\mbox{{\small -1}}}
\put(40,26){\mbox{{\small 3}}}
\put(35,26){\mbox{{\small 2}}}
\put(30,26){\mbox{{\small 1}}}

\put(10,15){\line(0,1){10}}
\put(15,15){\line(0,1){10}}
\put(20,25){\line(1,-1){10}}
\put(20,15){\line(1,1){4}}
\put(26,21){\line(1,1){4}}
\put(35,15){\line(0,1){10}}
\put(40,15){\line(0,1){10}}
\put(42,20){\mbox{$\cdots$}}

\linethickness{0.4mm}
\put(90,15){\line(0,1){3}}
\put(90,20){\line(0,1){5}}
\linethickness{0.2mm}
\put(88,21){\oval(4,4)[l]}
\put(88,19){\line(1,0){5}}
%\put(88,23){\line(1,0){5}}
\put(93,17){\oval(4,4)[tr]}
\put(93,25){\oval(4,4)[br]}
\put(95,26){\mbox{{\small 1}}}
\put(90,26){\mbox{{\small 0}}}

\put(100,15){\line(0,1){10}}

\put(107,20){\mbox{$\cdots$}}

\put(72,5){\mbox{$\tau_i$}}

\put(0,10){\line(1,-1){10}}
\put(0,0){\line(1,1){4}}
\put(6,6){\line(1,1){4}}
\put(45,10){\line(1,-1){10}}
\put(45,0){\line(1,1){4}}
\put(51,6){\line(1,1){4}}
\put(20,0){\line(0,1){10}}
\put(30,0){\line(0,1){10}}
\put(13,5){\mbox{$\cdots$}}
\put(35,5){\mbox{$\cdots$}}
\linethickness{0.4mm}
\put(90,0){\line(0,1){10}}
\linethickness{0.2mm}
\put(95,0){\line(0,1){10}}
\put(98,5){\mbox{$\cdots$}}
\put(102,10){\line(1,-1){10}}
\put(102,0){\line(1,1){4}}
\put(108,6){\line(1,1){4}}

\end{picture}
\caption{\label{generat} 
The graphical interpretation of the generators as symmetric tangles (on the left)
and as cylinder tangles (on the right)}

\end{figure}

$\ZB_n$ may be 
graphically interpreted (cf. figure \ref{generat}) as 
symmetric braids or cylinder braids: The symmetric picture
shows it as the group  of 
braids with $2n$ strands (numbered $-n,\ldots,-1,1,\ldots,n$) which are
fixed under a 180 degree rotation about the middle axis. 
In the cylinder picture one adds a single fixed line (indexed $0$)
on the left and
obtains $\ZB_n$ as the group of braids with $n$ strands that may 
surround this fixed line. 
The generators $\tau_i,i\geq0$ are mapped to the diagrams $X^{(G)}_i$ given in
figure \ref{generat}.

The braid group $\ZB_n$ has applications in the theory of knots in the
solid torus and in low dimensional physical systems with boundaries.
These applications motivate the search for tensor representations 
of $\ZB_n$ on n fold tensor product spaces $V^{\otimes n}$. 
It is natural to ask for extensions of the tensor representations of $\ZA_n$
given by quantum group R matrices. I. e. we are looking for
an endomorphism $F$ of $V$ such that the quantum braid matrix 
$B:=P(\pi\otimes\pi)R$ (P is the flip operator on $V\otimes V$ and 
$R$ is the
R matrix of the quantum group. $\pi$ is a representation on $V$.) 
fulfils
$F_1B_{1,2}F_1B_{1,2}=B_{1,2}F_1B_{1,2}F_1$ on $V^{\otimes n}$.
Subscripts indicate the spaces in which the matrices act.
$F$ is called the cylinder twist matrix.

 T. tom Dieck has found such extensions for the defining representations
 of the quantum groups of series $A,B,C$ in \cite{tD2}. 
 These solutions can be extended by   cabling to higher 
  representations as shown in \cite{TtDRHO}. 
  Naturally the question arises
 if these matrices come from an element in the quantum group.
 In \cite{TtDRHO} we have (from the point of view of universal
 operators) shown that this is the case
 for the quantum group of $sl_2$. 
 The present paper is a more detailed description of this
 result with calculations taking place in the quantum Weyl group.
We show that there is an element $t$ in the
 quantum Weyl group of $sl_2$  that gives rise to a cylinder twist 
 matrix $F=\pi(t)$ in every representation $\pi$. 
 This allows to calculate the representing matrices in all 
 dimensions.
  Taking the quantum Weyl group as starting algebra highlights  the Hopf
 algebraic content.

 The key observation behind our approach is the following:
The $F$ matrices of tom Dieck are all triangular with respect to the
counter diagonal. Hence taking out the quantum Weyl element as a factor 
one is left with a upper (or lower) triangular matrix. 
Such matrices may be representations matrices of an element
that is contained in one of the Borel sub algebras. Thus a simpler
ansatz may be used.

  \unitlength1mm

{\it Preliminaries:} Our notation for quantum groups is close to
\cite{KR}. The quantum group $\QG$ has relations
$[H,X]=2X, [H,Y]=-2Y,[X,Y]=(K^2-K^{-2})/(q^{1/2}-q^{-1/2})$
where $K:=q^{H/4},q:=e^h$. It follows that we have 
$KX=q^{1/2}XK,KY=q^{-1/2}YK$ and $Y^mH^n=(H+2m)^nY^m$. It is convenient to
use also $E:=KX, F:=K^{-1}Y$. One has 
$E^n=q^{-(n(n-1))/4}K^{n}X^n,F^n=q^{-(n(n-1))/4}K^{-n}Y^n$.
The coproduct is defined on generators by 
$\Delta(X):=X\otimes K+K^{-1}\otimes X ,
 \Delta(Y):=Y\otimes K+K^{-1}\otimes Y, \Delta(H)=H\otimes 1+1\otimes H$
 and the antipode by $S(H)=-H,S(X)=-q^{1/2}X,S(Y)=-q^{-1/2}Y$.

The associated quantum Weyl group $\QW$ is the Hopf algebra extension given
by an additional generator $w$ with relations
$wX=-q^{1/2}Yw,wY=-q^{-1/2}Xw,wH=-Hw$. It follows that
$wE=-q^{1/2}Fw,wF=-q^{-1/2}Ew$.
The Weyl element obeys $\epsilon(w)=1,w^2=v\epsilon$, 
where $v$ is the ribbon element and $\epsilon$ is $1$ in
odd dimensional and $-1$ in even dimensional irreducible
representations. 

The universal R matrix for both $\QG$ and $\QW$ is given by
\begin{equation}
R=q^{H\otimes H/4}\sum_{n=0}^\infty\frac{(1-q^{-1})^n}{[n]!}
q^{n(n-1)/4}E^n\otimes F^n
\end{equation}
Here we have used the usual quantum factorial defined from the quantum
number $[n]:=(q^{n/2}-q^{-n/2})/(q^{1/2}-q^{-1/2})$.
We have to calculate the antipode of $w$ because there are 
different results stated in the literature \cite{KR},\cite{LS}.
We write $R=\sum_k\Ralpha_k\otimes\Rbeta_k$ and 
introduce $c_n:=\frac{(1-q^{-1})^n}{[n]!}q^{n(n-1)/4},
u:=\sum_k S(\Rbeta_k)\Ralpha_k$. 
\begin{eqnarray}
  1&=&\epsilon(w^{-1})=m(S\otimes \id)\Delta(w^{-1})
  =m(S\otimes \id)((w^{-1}\otimes w^{-1})R)\nonumber\\
  &=&m(S\otimes \id)(R_{2,1}(w^{-1}\otimes w^{-1}))
  =\sum_k m(S\otimes\id)(\Rbeta_kw^{-1}\otimes\Ralpha_kw^{-1})\nonumber\\
  &=&\sum_k S(w)^{-1}S(\Rbeta_k)\Ralpha_kw^{-1}=S(w^{-1})uw^{-1}\nonumber\\
  S(w)&=&uw^{-1}
\end{eqnarray}
Essential for our calculations is the formula for the 
coproduct of the Weyl element 
$\Delta(w)=R^{-1}(w\otimes w)$ and a simple implication
$(w\otimes w)R=R_{2,1}(w\otimes w)$.

\section{Construction of the Cylinder twist}

In this section we construct a solution of the cylinder braid
equation
\begin{equation}
R_{2,1}t_2Rt_1=t_1R_{2,1}t_2R			  \label{vierzopf}
\end{equation}
As motivated in the introduction we use the ansatz $t=wz$.
\begin{lemma} \label{zlem} Equation (\ref{vierzopf}) holds with $t=wz$ if
\begin{equation} \Delta(z)=z_2w_2R_{2,1}w_2^{-1}z_1\label{zdelta}
\end{equation}
\end{lemma}
\begin{bew}	We express $R_{2,1}$ using $w$.
\begin{eqnarray*}
\lefteqn{(\ref{vierzopf})
\Leftrightarrow w_1w_2Rw_1^{-1}w_2^{-1}w_2z_2Rw_1z_1=
  w_1z_1w_1w_2Rw_1^{-1}w_2^{-1}w_2z_2R}\\
&\Leftrightarrow&Rw_1^{-1}z_2Rw_1z_1=z_1w_1Rw_1^{-1}z_2R\\
&\Leftrightarrow&Rw_1^{-1}z_2 w_1w_2R_{2,1}w_1^{-1}w_2^{-1} w_1z_1=
z_1w_1Rw_1^{-1}z_2R\\
&\Leftrightarrow&Rz_2 w_2R_{2,1}w_2^{-1} z_1=z_1w_1Rw_1^{-1}z_2R\\
&\Leftrightarrow& R\Delta(z)=\Delta'(z)R
\end{eqnarray*}
The last line holds because $R$ intertwines between the coproduct and the
opposite coproduct.
\end{bew}

If one had made the ansatz $t=w^{-1}z$ the condition would be
$\Delta(z)=z_2w_2^{-1}R_{2,1}w_2z_1$ and $t=zw$ would lead
to $\Delta(z)=z_1w_1R_{2,1}w^{-1}_1z_2$.

Note  that $t':=w^{-2}t=w^{-1}z$ is another solution of (\ref{vierzopf})
because $w^2$ is central.

\begin{bem}
If $z$ is a solution of (\ref{zdelta}) then so is 
$\widetilde{z}:=K^\alpha z K^\alpha$ where $\alpha$ 
is an arbitrary number. The computation is straightforward
from the fact that $R$ and $K^\alpha\otimes K^\alpha$ commute.
\end{bem}

\begin{bem}
If $z$ is a solution of (\ref{zdelta}) then so is $S(u)^{-1}zu$. 
To prove this we first note that $uw=S^2(w)u=S(uw^{-1})u=S(w^{-1})S(u)u=
wu^{-1}S(u)u=Cwu^{-1}=wS(u)$ where $C:=uS(u)=S(u)u$ is 
Drinfeld's Casimir operator. A consequence is
$Cu^{-1}w=wu$. We first investigate the behaviour of $uzu$:
\begin{eqnarray*}
\lefteqn{(1\otimes uzuw)R_{2,1}(uzu\otimes w^{-1})
= (uu^{-1}\otimes uzCwu^{-1})R_{2,1}(uzu\otimes uu^{-1}w^{-1})}\\
&=&	(u\otimes u)(1\otimes Czw)R_{2,1}(zu\otimes C^{-1}w^{-1}u)
= (u\otimes u)(1\otimes zw)R_{2,1}(z\otimes w^{-1})(u\otimes u)\\
&=&(u\otimes u)\Delta(z)(u\otimes u)=
R_{2,1}R\Delta(u)\Delta(z)\Delta(u)R_{2,1}R=\Delta(uzu)(R_{2,1}R)^2
\end{eqnarray*}
Here we have used that $u\otimes u$ and $R$ commute.
\begin{eqnarray*}
& &(1\otimes C^{-1}uzuw)R_{2,1}(C^{-1}uzu\otimes w^{-1})
=(C^{-1}\otimes C^{-1})(1\otimes uzuw)R_{2,1}(uzu\otimes w^{-1})\\
&=&
(C^{-1}\otimes C^{-1})(R_{2,1}R)^2\Delta(uzu)=\Delta(C^{-1})\Delta(uzu)=
\Delta(C^{-1}uzu)
\end{eqnarray*}
Now, $C^{-1}uzu=S(u)^{-1}zu$ and the claim is shown.
\end{bem}

\begin{bem}
The element $t$ gives not only rise to representations of $\ZB_n$
but also of the braid group of the affine series $C_n^{(1)}$. 
Tensor representations of this braid group need another element $\overline{F}$
such that $(1\otimes\overline{F})B(1\otimes\overline{F})B=
B(1\otimes\overline{F})B(1\otimes\overline{F})$ in $V\otimes V$.
If $\overline{F}=\pi(\overline{t})$ this is equivalent to
$R_{2,1}\overline{t}_1R\overline{t}_2=\overline{t}_2R_{2,1}\overline{t}_1R$
which is (by permuting the tensor factors)
equivalent to
$R\overline{t}_2R_{2,1}\overline{t}_1=\overline{t}_1R\overline{t}_2R_{2,1}$.
Now, assume that $t$ is any solution of (\ref{vierzopf})
and multiply (\ref{vierzopf}) from the left with $(w_1w_2)^{-1}$
and from the right with $w_1w_2$. This shows that $\overline{t}=w^{-1}tw$
provides a solution of the above equation.
\end{bem}

The element $w_2R_{2,1}w_2^{-1}$ that occurs in 
(\ref{zdelta}) is explicitly:
\begin{equation}
w_2R_{2,1}w_2^{-1}=q^{-H\otimes H/4}\sum_{n=0}^\infty
\frac{(1-q^{-1})^n}{[n]!}
q^{n(n-1)/4}(-q^{1/2})^n F^n\otimes F^n		 \label{rdach}
\end{equation}

This shows that we may assume that $z$ is an element of 
the Borel sub-algebra generated by $H,Y$.
In order to reproduce the first factor in (\ref{rdach})
from a coproduct it seems to be adequate to make the further
factorisation $z=q^{-H^2/8}\widehat{z}$. Note that the factor
$q^{-H^2/8}$ already occurred in \cite{KR} in the connection with
Lusztig's automorphisms.
\begin{equation}
\Delta(z)=\left(q^{-H^2/8}\otimes q^{-H^2/8}\right)
q^{-H\otimes H/4}\Delta(\widehat{z})  \label{leq1}
\end{equation}
The right-hand side of (\ref{zdelta}) becomes
\begin{eqnarray}
\lefteqn{ (1\otimes q^{-H^2/8})(1\otimes\widehat{z})
q^{-H\otimes H/4} 
\sum_{n=0}^\infty
\frac{(1-q^{-1})^n}{[n]!}
q^{n(n-1)/4}(-q^{1/2})^n }\nonumber\\
&&(F^n\otimes F^n)		
 (1\otimes q^{-H^2/8})(1\otimes\widehat{z})	\nonumber \\
&=&(q^{-H^2/8}\otimes q^{-H^2/8})(1\otimes\widehat{z})
q^{-H\otimes H/4} \nonumber\\
&&\sum_{n=0}^\infty
\frac{(1-q^{-1})^n}{[n]!}
q^{n(n-1)/4}(-q^{1/2})^n q^{-n^2/2} ((q^{-Hn/2}F^n)\otimes F^n)		
 (\widehat{z}\otimes 1)						  \label{req1}
\end{eqnarray}
Here we have used $F^nq^{-H^2/8}=q^{-(H+2n)^2/8}F^n$.

Cancelling the leftmost factors in $(\ref{leq1})=(\ref{req1})$
and introducing the shortcut $B_n:=	\frac{(1-q^{-1})^n}{[n]!}
q^{n(n-1)/4}(-q^{1/2})^n q^{-n^2/2}$ we arrive at the following
equation which has to be solved
\begin{equation}
\Delta(\widehat{z}) = q^{H\otimes H/4} (1\otimes\widehat{z})
q^{-H\otimes H/4} 
\sum_{n=0}^\infty
B_n ((q^{-Hn/2}F^n)\otimes F^n)
(\widehat{z}\otimes 1) \label{glg}		
\end{equation}

It seems to be difficult to go on with a general ansatz 
$\sum_{i,j}c_{i,j}K^iF^j$ for $\widehat{z}$.
Experiments show that the following ansatz works:
\begin{equation}
\widehat{z}=\sum_{m=0}^\infty\beta_m q^{\alpha Hm}Y^m \label{zhansatz}
\end{equation}
Here $\alpha,\beta_m$ are coefficients which are yet to be determined.
It is important not to use $F^m$ in (\ref{zhansatz}) because 
this would cause the coproduct to produce unbalanced $K$ factors.

\begin{eqnarray}
\Delta(\widehat{z})&=&\sum_m\beta_m(q^{\alpha Hm}\otimes
q^{\alpha Hm})(Y\otimes K+K^{-1}\otimes Y)^m\nonumber\\
&=&\sum_m\sum_{i=0}^m\beta_m\frac{[m]!}{[i]![m-i]!}q^{i(m-i)/2}
(q^{\alpha Hm}\otimes q^{\alpha Hm}) 
(Y\otimes K)^i(K^{-1}\otimes Y)^{m-i}\nonumber\\
&=&	\sum_m\sum_{i=0}^m\beta_m\frac{[m]!}{[i]![m-i]!}
(q^{\alpha Hm}\otimes q^{\alpha Hm})
q^{(i-m)H/4}Y^i \otimes q^{iH/4}Y^{m-i}	\label{RL}
\end{eqnarray}
The right-hand side of (\ref{glg}) becomes
\begin{eqnarray}
\lefteqn{
\sum_{s,t}\beta_s\beta_t\sum_n B_nq^{H\otimes H/4}
(1\otimes q^{\alpha Hs}Y^s)q^{-H\otimes H/4} }\nonumber\\
&&(q^{-Hn/2}(K^{-1}Y)^n\otimes(K^{-1}Y)^n)(q^{\alpha Ht}Y^t\otimes1)
 \nonumber\\
&=&	\sum_{s,t}\beta_s\beta_t\sum_n B_nq^{-n(n-1)/2}
(q^{-Hs/2}\otimes q^{\alpha Hs}Y^s)
(q^{-Hn/2}K^{-n}Y^nq^{\alpha Ht}Y^t\otimes K^{-n}Y^n)
\nonumber\\
&=&	\sum_{s,t}\beta_s\beta_t\sum_n B_nq^{-n(n-1)/2}q^{-sn/2+2n\alpha t}
\nonumber\\&&
(q^{-Hs/2}q^{-Hn/2}K^{-n}q^{\alpha Ht}Y^{n+t}\otimes q^{\alpha Hs}K^{-n}Y^{s+n})
\nonumber\\
&=&	\sum_{s,t}\beta_s\beta_t\sum_n B_nq^{-n(n-1)/2-sn/2+2n\alpha t}
\nonumber\\&&
(q^{H(-s-n-n/2+2\alpha t)/2}Y^{n+t}\otimes q^{H(2\alpha s-n/2)/2}Y^{s+n})
\label{LL}
\end{eqnarray}

Since $H^iY^j$ is a basis of the Borel sub-algebra we can make
a term by term comparison of the coefficients of $Y^a\otimes Y^b$.
We start by investigating the first few terms:
\begin{equation}\begin{array}{lccccccccr}
a&b  & i & m & \mbox{Coeff. (\ref{RL})}
     & n & t & s & \mbox{Coeff. (\ref{LL})}	& \mbox{Conclusion}\\
0&0 & 0 & 0 & \beta_0 & 0& 0& 0 & \beta_0^2B_0 & \beta_0=1\\
0&1 & 0&1 & \beta_1q^{H(\alpha-1/4)}\otimes q^{\alpha H} &
0&0&1 & \beta_1q^{-H/2}\otimes q^{\alpha H} & \alpha=-1/4 \\
1&0 & 1 & 1 & \beta_1(q^{-H/4}\otimes 1) & 0&1&0 &
\beta_1q^{-H/4}\otimes 1 & -
\end{array}\end{equation}
Now that we have determined $\alpha$ and $\beta_0$ we can consider the general
case. The coefficient of $Y^a\otimes Y^b$ in (\ref{RL}) is
\begin{eqnarray*}
\lefteqn{\beta_{a+b}\frac{[a+b]!}{[a]![b]!}(q^{-H(a+b)/4}\otimes q^{-H(a+b)/4})
(q^{-bH/4}\otimes q^{aH/4})}\\
&&=	 \beta_{a+b}\frac{[a+b]!}{[a]![b]!}(q^{-H(a/4+b/2)}\otimes q^{-Hb/4})
\end{eqnarray*}
In (\ref{LL}) we set $t=a-n\geq0,s=b-n\geq0$ and obtain the coefficient:
\begin{eqnarray*}
\lefteqn{\sum_n B_n\beta_{a-n}\beta_{b-n}
        q^{-n(n-1)/2-(b-n)n/2-n(a-n)/2}}\\
 && q^{H(n-b-n-n/2-(a-n)/2)/2}\otimes q^{H(-(b-n)/2-n/2)/2}\\
&&=	\sum_n B_n\beta_{a-n}\beta_{b-n}q^{n^2/2-n(a+b-1)/2}
  q^{H(-b-a/2)/2}\otimes q^{-Hb/4}
\end{eqnarray*}

The terms involving $H$ are equal. We are left with 
\begin{equation}   \label{bed}
  \beta_{a+b}\frac{[a+b]!}{[a]![b]!}=	
  \sum_n B_n\beta_{a-n}\beta_{b-n}q^{n^2/2-n(a+b-1)/2}
\end{equation}
Since we have $s=b-n,t=a-n\geq0$ the sum over $n$ runs only from $0$ to 
${\rm min}(a,b)$. We set $b=1$ and obtain
\begin{eqnarray*}
\beta_{a+1}[a+1]&=&\beta_a\beta_1+B_1\beta_{a-1}q^{1/2-a/2}	\\
&=&	\beta_a\beta_1+\beta_{a-1}(q^{-1}-1)q^{(1-a)/2}
\end{eqnarray*}

It remains to show that (\ref{bed}) holds for $b>1$. 
To this end we first simplify the recursion formula by 
defining $\beta'_a:=\beta_a[a]!$:
\begin{eqnarray}
\beta'_{a+1}&=&\beta'_1\beta'_a+\beta'_{a-1}(q^{-a}-1)
\end{eqnarray}
We first reformulate  (\ref{bed})
it in terms of $\beta'_a$:
\begin{eqnarray}
\beta'_{a+b}&=&\sum_{n=0}^{{\rm min}(a,b)} B_n
\frac{[a]!}{[a-n]!}\frac{[b]!}{[b-n]!}\beta'_{a-n}\beta'_{b-n}
q^{n^2/2-n(a+b-1)/2}\nonumber\\
&=&	\sum_{n=0}^{{\rm min}(a,b)} 
\left[{a \atop n}\right] \left[{b\atop n}\right][n]!
\beta'_{a-n}\beta'_{b-n}
q^{n^2/2-n(a+b-1)/2}q^{n(n-1)/4}(q^{-1}-1)^nq^{n/2-n^2/2}\nonumber\\
&=&	\sum_{n=0}^{{\rm min}(a,b)} B^{a,b}_n 	\beta'_{a-n}\beta'_{b-n}
\label{bform}\\
B^{a,b}_n&:=&  \left[{a \atop n}\right] \left[{b\atop n}\right][n]!
q^{-n(a+b)/2}q^{3n/4+n^2/4}(q^{-1}-1)^n
\end{eqnarray}
We are done if we can show that	 the right hand side of (\ref{bform})
is actually only a function of $a+b$. This will follow from the
fact that the substitutions $a\rightarrow a+1$ and $b\rightarrow b+1$
have the same effect.
We first concentrate on the substitution $a\rightarrow a+1$ 
and calculate the coefficient $B^{a+1,b}_n$ using the formula
for the q-binomial coefficients:
\begin{eqnarray*}
\left[{a+1\atop n}\right]&=&q^{-n/2}\left[{a\atop n}\right]+
q^{(a+1-n)/2}\left[{a\atop n-1}\right]\\
B^{a+1,b}_n
&=&q^{-n}B^{a,b}_n+q^{(a+1-n)/2}\left[{a\atop n-1}\right]
\left[{b\atop n}\right][n]!q^{-n(a+1+b)/2}q^{3n/4+n^2/4}(q^{-1}-1)^n\\
&=&q^{-n}B^{a,b}_n+q^{(a+1-n)/2}\left[{a\atop n-1}\right]
\left[{b\atop n-1}\right][n-1]![b-n+1]\\&&
q^{-n(a+1+b)/2+3n/4+n^2/4}(q^{-1}-1)^n\\
&=&q^{-n}B^{a,b}_n+B^{a,b}_{n-1}[b-n+1](q^{-1}-1)q^{-n/2}
q^{(a+1-n)/2}q^{3/4+n/2-1/4-(a+b)/2}\\
&=&q^{-n}B^{a,b}_n+B^{a,b}_{n-1}q^{(1-n)/2}q^{(1-b)/2}(q^{-1}-1)
[b-n+1]\\
&=&q^{-n}B^{a,b}_n+q^{(1-n-b)/2}(q^{(n-b-1)/2}-q^{(b-n+1)/2})B^{a,b}_{n-1}\\
&=&q^{-n}(B^{a,b}_n+(q^{n-b}-q)B^{a,b}_{n-1}) \\
B^{a,b+1}_n&=& q^{-n}(B^{a,b}_n+(q^{n-a}-q)B^{a,b}_{n-1})
\end{eqnarray*}
We now consider the righthand side of (\ref{bform}).
It is convenient to set $B^{a,b}_n=\beta'_n=0$ for negative $n$.
Doing this we don't have to care about summation ranges
and can freely shift the summation variable as we do in the third
step and in the first summand in the fifth step
of the following calculation:
\begin{eqnarray*}
\lefteqn{\sum_n B^{a+1,b}_n\beta'_{a+1-n}\beta'_{b-n}=}\\
&=&\sum_n q^{-n}(B^{a,b}_n+(q^{n-b}-q)B^{a,b}_{n-1})
\beta'_{a-n+1}\beta'_{b-n}\\
&=&\sum_n q^{-n-1}(B^{a,b}_{n+1}+(q^{n+1-b}-q)B^{a,b}_n)
\beta'_{a-n}\beta'_{b-n-1} \\
&=&	\sum_nq^{-n-1}B^{a,b}_{n+1}\beta'_{a-n}\beta'_{b-n-1}
+\sum_nq^{-n}B^{a,b}_n\beta'_{a-n}(q^{n-b}-1)\beta'_{b-n-1}\\
&=&	\sum_nq^{-n}B^{a,b}_{n}\beta'_{a-n+1}\beta'_{b-n}
+\sum_nq^{-n}B^{a,b}_n\beta'_{a-n}(\beta'_{b-n+1}-\beta'_1\beta'_{b-n})\\
&=&	\sum_nq^{-n}B^{a,b}_{n}\beta'_{a-n+1}\beta'_{b-n}
+ \sum_nq^{-n}B^{a,b}_n\beta'_{a-n}\beta'_{b-n+1}
-\sum_nq^{-n}B^{a,b}_n\beta'_{a-n}\beta'_1\beta'_{b-n}
\end{eqnarray*}
The calculation of $ \sum_n B^{a,b+1}_n\beta'_{a-n}\beta'_{b+1-n}$
would give the same terms: Exchanging $a$ and $b$ interchanges
the
first two summands and leaves the third invariant.
The proof is complete.

\begin{satz}
A solution of $R_{2,1}t_2Rt_1=t_1R_{2,1}t_2R$ is given by 
\begin{equation}
t=wq^{-H^2/8}\sum_{m=0}^\infty\beta_mq^{-Hm/4}Y^m\end{equation}
where $\beta_0=1$, $\beta_1$ is arbitrary and 
\begin{equation}
\beta_{a+1}=(\beta_a\beta_1+\beta_{a-1}(q^{-1}-1)q^{(1-a)/2})/[a+1]
\end{equation}\end{satz}

From this expressions $t$ can be calculated in all irreducible
representations of $\QW$. We give the matrices in the $2$, $3$ and $4$
dimensional representations. We use the standard basis in
the order of decreasing weights.
\begin{eqnarray}
\left(\begin{array}{cc}
-\beta_1q^{-1/2} & -q^{-3/4} \\ q^{-1/4} & 0
\end{array}\right)	\\
\left(\begin{array}{ccc}
(1-q+q\beta_1^2)/q^2 & q^{-7/4}(q+1)^{1/2}\beta_1 & q^{-2} \\
-q^{-5/4}\beta_1(q+1)^{1/2} & -q^{-1} & 0 \\
q^{-1} & 0 & 0
\end{array}\right)	 \\
\left(\begin{array}{cccc}
-q^{-7/2}\beta_1(1+q-2q^2+q^2\beta_1^2) &
q^{-15/4}\gamma(q-1-q\beta_1^2) &
 -q^{-7/2}\gamma\beta_1& -q^{-15/4} \\
 q^{-13/4}\gamma(1-q+q\beta_1^2)&q^{-5/2}(1+q)\beta_1 & q^{-9/4} & 0 \\
-q^{-5/2}\gamma\beta_1 & -q^{-7/4} & 0 & 0
\\ q^{-9/4} & 0 & 0 & 0
\end{array}\right) \\
\gamma:= (1+q+q^2)^{1/2}
\end{eqnarray}
It can easily be checked that these matrices indeed fulfil
(\ref{vierzopf}). 

It should be noted that the proposition leads to a second infinite 
series of tensor representations  of  the braid group $\ZB_n$ because
there is a second series of irreducible representations 
of the quantum Weyl group $\QW$. These representations are not 
irreducible as representations of $\QG$.

We calculate the inverse of $\widehat{z}$.
 \begin{lemma}
 \begin{eqnarray}
 \widehat{z}^{-1}&=&\sum_{m=0}^\infty\alpha_m q^{-Hm/4}Y^m\\
 \alpha_0&=&1\\
 \alpha_a&=&-\sum_{m=0}^{a-1}\alpha_{a-1-m}\beta_mq^{-m(a-1-m)/2}
\end{eqnarray}
 \end{lemma}
\begin{bew}
\begin{eqnarray*}
1=\widehat{z}\widehat{z}^{-1}&=&
\sum_{m,n}\alpha_n\beta_m q^{-Hm/4}Y^mq^{-Hn/4}Y^n\\
&=&	\sum_{m,n}\alpha_n\beta_m q^{-Hm/4}q^{-(H+2m)n/4}Y^{m+n}\\
&=&	\sum_{m,n}\alpha_n\beta_m q^{-nm/2}q^{-H(m+n)/4}Y^{m+n}
\end{eqnarray*}
Looking at the term with $Y^0$ one obtains $\alpha_0=1$. 
We now consider the terms with $q^{-Ha/4}Y^a,a\geq1$. We substitute $n=a-m$.
Since $n\geq0$ we obtain a restriction for the $m$ summation
$0\leq m\leq a$
The coefficient that should vanish is 
\[\sum_{m=0}^a\alpha_{a-m}\beta_mq^{-m(a-m)/2}\]
Isolation of $\alpha_a$ yields the formula given in the lemma.
\end{bew}

\section{Properties of the cylinder twist}

In this section we try to analyse the the algebraic properties of $t$ and 
try to fit them in a broader framework.
 
 In \cite{rhobcat} we were led by categorial considerations to the following
 axioms: A restricted Coxeter-B braided Hopf algebra 
 is a ribbon Hopf algebra with
  an element $\overline{v}\in H$
 such that
 \begin{eqnarray}
  R_{2,1}\overline{v}_2R\overline{v}_1&=&
 \overline{v}_1R_{2,1}\overline{v}_2R	\label{vierv}  \\
   \epsilon(\overline{v})&=&1\label{hd1}\label{vepsilon}\\
 \Delta(\overline{v})&=&R^{-1}\overline{v}_2R\overline{v}_1\label{vdelta} \\
   S(\overline{v})&=&v^2\overline{v}^{-1}		\label{hd3}
 \end{eqnarray}

Now, we consider the properties of $t$.
 \begin{satz}\begin{eqnarray}
 R_{2,1}t_2Rt_1&=&t_1R_{2,1}t_2R	  \\
 \Delta(t)&=&R^{-1}t_2Rt_1\\
 \epsilon(t)&=&1
 \end{eqnarray}
 \end{satz}
 \begin{bew}
 The first equation has already been proven. 
\begin{eqnarray}
\Delta(t)&=&\Delta(w)\Delta(z)=
R^{-1}w_1w_2 z_2w_2R_{2,1}w_2^{-1}z_1\nonumber\\
&=&R^{-1}w_2z_2w_1w_2R_{2,1}w_2^{-1}w_1^{-1}w_1z_1
=R^{-1}w_2z_2Rw_1z_1=R^{-1}t_2Rt_1 \nonumber
\end{eqnarray}
The third equation is trivial.
\end{bew}

Since (\ref{hd3}) follows from the remaining axioms we conclude that
the quantum Weyl group of $sl_2$ is a restricted Coxeter-B braided
Hopf algebra.

\section{Outlook}

A natural further challenge is to find a cylinder twist element $t$
in the quantum Weyl groups associated to other Lie algebras.
The Weyl element  $w$ has a natural generalisation as the longest element
$w_0$ in the quantum Weyl group. However, generalising our
construction of $z$ would  require to evaluate products 
of root vectors which is a highly non-trivial task.
We would need a sort of quantum double construction of $z$. 
Using the fact that $w_0$ maps a positive root vector to a multiple of a 
negative root vector and hence interchanges $H$ and $H^\ast$
in the quantum double construction one may write down 
a version of (\ref{zdelta}). 
To be more precise, consider the quantum double realised
as in \cite{majid} on $H^\ast\otimes H$. With dual bases $\{f^a\},\{e^a\}$
the $R$ matrix is $R=\sum_a f^a\otimes1\otimes1\otimes e^a$. 
We assume that there are scalars $\lambda_a$ such that
$w(f^a\otimes 1)w^{-1}=\lambda_a(1\otimes e^a)$. This implies
$w(1\otimes e^a)w^{-1}=\lambda_a^{-1}(f^a\otimes1)$ and 
 $(w\otimes w)R=R_{2,1}(w\otimes w)$.
Furthermore, we assume $\Delta(w)=R^{-1}(w\otimes w)$.
Then  the ansatz $z=\sum_a\beta_a(1\otimes e^a)$ turns 
equation (\ref{zdelta}) into
\[
 \sum_a\beta_a(1\otimes e^a_1\otimes1\otimes e^a_2)=
 \sum_{a,b,n}\lambda_n\beta_a\beta_b(1\otimes e^ne^a\otimes1\otimes e^be^n)
\]
Unfortunately, no natural solution suggests itself.

Furthermore, one should clarify possible connections with other
occurrences of (\ref{vierzopf}), especially Majid's theory
of braided Lie algebras \cite{majid}.

 %\section{Towards a generalisation to the general case}

 %Kirillov and Reshetikhin introduce the general quantum Weyl group
 %associated to the quantum group of an arbitrary semsimple Lie algebra.
 %They use two versions of the Weyl elements 
 %$\check{w}_i:=w_iq_i^{H_i^2/8}$. The Lusztig automorphisms are then
 %$T_i(a)=\check{w}_i^{-1}a\check{w}_i$. 
 %The analogon of the longest element in the Weyl group is 
 %denoted by $\check{w}_0$. It obeys:
 %\begin{equation}
 %R=\exp\left(h/2\sum_{i,j}^nB^{-1}_{i,j} H_i\otimes H_j\right)
 %(\check{w}_0\otimes\check{w}_0)\Delta(\check{w}_0)^{-1} \label{rr0}
 %\end{equation}
 %Since the exponential term in (\ref{rr0}) is symmetric
 %one can deduce that $R(\check{w}_0\otimes\check{w}_0)=
 %(\check{w}_0\otimes\check{w}_0)R_{2,1}$ and thus
 %lemma \ref{zlem} can be applied with $w=\check{w}_0$.
 %Moreover, the remarks following lemma \ref{zlem} are 
 %valid as well.

\end{document}